\documentclass[leqno,12pt]{article}

\usepackage{amstext    }
\usepackage{amsthm    }
\usepackage{a4}
\usepackage[mathscr]{eucal}
\usepackage{mathrsfs}

\usepackage{amsmath}
\usepackage{amssymb}
\usepackage{amscd}

\numberwithin{equation}{section}

\newtheorem{theorem}{Theorem}[section]
\newtheorem{definition}[theorem]{Definition}
\newtheorem{proposition}[theorem]{Proposition}
\newtheorem{corollary}[theorem]{Corollary}
\newtheorem{lemma}[theorem]{Lemma}

\newtheorem{example}[theorem]{Example}

\newtheorem{problem}[theorem]{Problem}

\newcommand{\cali}[1]{\mathscr{#1}}

\newcommand{\GL}{{\rm GL}}

\newcommand{\SL}{{\rm SL}}

\newcommand{\Aut}{{\rm Aut}}
\newcommand{\Alb}{{\rm Alb}}

\newcommand{\volume}{{\rm vol}}

\renewcommand{\Im}{\mathop{\mathrm{Im}}}

\newcommand{\dist}{\mathop{\mathrm{dist}}\nolimits}

\renewcommand{\GL}{{\rm GL}}

\newcommand{\id}{{\rm id}}

\newcommand{\Kc}{\cali{K}}

\newcommand{\C}{\mathbb{C}}

\newcommand{\N}{\mathbb{N}}
\newcommand{\Z}{\mathbb{Z}}
\newcommand{\R}{\mathbb{R}}
\newcommand{\Q}{\mathbb{Q}}

\renewcommand\P{\mathbb{P}}

\newcommand{\A}{\mathbb{A}}
\newcommand{\F}{\mathbb{F}}

\title{Tits alternative for automorphism groups of compact K\"ahler manifolds}

\author{Tien-Cuong Dinh}

\begin{document}

\maketitle

\begin{abstract}
We survey some properties of the automorphism groups of compact K\"ahler manifolds. In particular, we present recent results by Keum, Oguiso and Zhang on the structure of these groups from the Tits alternative point of view. Several other related results will be also discussed.  
\end{abstract}

\noindent
{\bf Classification AMS 2010 :} 14J50, 14E07, 32M05, 32H50, 32Q15.


\section{Introduction} \label{intro}

In this paper, we will present a proof of the Tits type alternative which confirms a conjecture by Keum-Oguiso-Zhang \cite{KOZ}. This nice result was proved in \cite{KOZ} for some particular cases and in Zhang's work
 \cite{Zhang3} for the general case. It 
is intimately connected to  developments in complex dynamics of several variables. So we will survey some results from holomorphic dynamics which are related to the Tits type alternative and to its proof.

Several statements in the paper are purely algebraic and one can ask if they can be obtained by purely algebraic methods and also  if they can be extended to automorphisms of algebraic manifolds over a finite field. The question is far from the author's competence. So only analytical tools will be discussed here and we refer the reader to a recent paper by Esnault-Srinivas \cite{ES} which is a  first step to study similar questions for finite fields.

Let $X$ be a compact K\"ahler manifold, e.g. a complex projective manifold, of dimension $k$. 
Denote by $\Aut(X)$ the group of all holomorphic automorphisms of $X$. 
Following a result by Bochner-Montgomery,
$\Aut(X)$ is a complex Lie group of finite dimension \cite{BM}, see also \cite{Akhiezer,BM,SW} and the references therein for results on upper bounds of the dimension. 

The group $\Aut(X)$ may have an infinite number of connected components. Let $\Aut_0(X)$ denote the connected component of the identity. Elements in this subgroup are those induced by global holomorphic vector fields on $X$. They are almost characterized by the property that the associated actions on cohomology preserve a K\"ahler class. More precisely, Fujiki and Lieberman proved that the group of all automorphisms preserving a given K\"ahler class is a finite union of connected components of $\Aut(X)$, see \cite{Fujiki,Lieberman} and Theorem \ref{t:FL1} below.

Quite recently, ideas from complex dynamics allowed to study the "discrete part" of $\Aut(X)$. Inspired by results in \cite{DS1}, a Tits type alternative was conjectured and proved in some cases by Keum-Oguiso-Zhang in \cite{KOZ}, see also \cite{Oguiso,Zhang1}.  The conjecture was fully obtained by De-Qi Zhang in \cite{Zhang3}. His theorem corresponds to the second assertion of the following theorem, see \cite[Lemma 2.11]{Zhang3}. The first assertion was obtained by Cantat in \cite{Cantat} for finite type groups and Campana-Wang-Zhang in \cite{CWZ} for the general case.

\begin{theorem} \label{t:Zhang} 
Let $X$ be a compact K\"ahler manifold of dimension $k$ and of Kodaira dimension $\kappa_X$. 
Define $\kappa:=\max(\kappa_X,0)$ if $\kappa_X<k$ and $\kappa:=k-1$ otherwise. 
Let $G$ be a group of holomorphic automorphisms of 
$X$ which does not contain any free non-abelian subgroup. Then $G$ admits a finite index subgroup $G'$ satisfying the following properties:
\begin{enumerate}
\item $G'$ is solvable; in other words, $\Aut(X)$ satisfies the Tits alternative, see also Theorem \ref{t:Tits} below;
\item The set $N'$ of zero entropy elements of $G'$  is a normal subgroup of $G'$ and $G'/N'$ is a free abelian group of rank at most  $k-\kappa-1$. 
\end{enumerate}
\end{theorem}

The entropy of an automorphism was originally introduced as a dynamical invariant. However, thanks to results by Gromov and Yomdin \cite{Gromov2, Yomdin}, it can be also considered as an algebraic invariant. The notion, its properties and its relations with the dynamical degrees of automorphism will be presented in Section \ref{s:entropy}. 
Note that in Theorem \ref{t:Zhang}, the fact that $N'$ is a subgroup of $G'$ is not obvious.  

A main tool in the proof of Theorem \ref{t:Zhang} is a mixed version of the classical Hodge-Riemann theorem. It will be discussed in Section \ref{s:HR}. 
In Section \ref{s:fibration}, we will survey some results on meromorphic fibrations which are preserved by an automorphism. These results will be applied to the Iitaka's fibration of $X$ and are used to obtain the rank estimate in Theorem \ref{t:Zhang}. 
Finally, the proof of Theorem \ref{t:Zhang} will be given in Section \ref{s:Tits}. 

\begin{problem}\rm 
Can we always choose $G'$ so that its derived length is bounded by a constant which depends only on the dimension $k$ ?
\end{problem}

Using the Hodge-Riemann and the Lie-Kolchin theorems, we can show that the answer is affirmative for $k=2$.

We deduce from Theorem \ref{t:Zhang} the following consequence, see \cite{DS1,Zhang3}. We say that a group $G$ of automorphisms of $X$ has {\it positive entropy} if all elements of $g$, except the identity, have positive entropy. 

\begin{corollary} \label{c:commute}
Let $G$ be a group of automorphisms of $X$. Assume that $G$ is abelian and has positive entropy. Then, $G$ is a free abelian group of rank at most  $k-\kappa-1$. 
\end{corollary}

The rank estimates in Theorem \ref{t:Zhang} and Corollary \ref{c:commute} are optimal as shown in  the following example. 

\begin{example} \rm \label{ex_torus}
Consider the natural action of $\SL(k,\Z)$ on the complex torus $\C^k/(\Z^k+i\Z^k)$. The Kodaira dimension of the torus is zero. By a theorem of  Prasad-Raghunathan \cite{PR}, the action on the right of 
$\SL(k,\Z)\setminus \SL(k,\R)$
of the group of diagonal matrices in $\SL(k,\R)$ admits compact orbits. These compact orbits can be identified to quotients of $\R^{k-1}$ by subgroups of $\SL(k,\Z)$. We deduce that $\SL(k,\Z)$ admits free abelian subgroups of rank $k-1$ which are diagonalizable. The elements of these subgroups, except the identity, admit eigenvalues with modulus larger than 1.  
They can be identified to rank $k-1$ free abelian groups of automorphisms on the considered complex torus. It is then not difficult to check that such groups have positive entropy.

The existence of  free abelian subgroups of $\SL(k,\Z)$ can be also obtained in the following way (proof communicated to the author by B. Gross). 
Let $\A$ be the ring of algebraic integers of a totally real number field $\F$ of degree $k$. 
We have  $\F\otimes_\Q \R=\R^k$.
Let $\A^\times$ denote the group of units in $\A$. By Dirichlet's unit theorem, $\A^\times$ admits a free abelian subgroup $G$ of rank $k-1$, see \cite{Lang}. The multiplication in $\F$  induces an injective representation $\rho:G\to \GL(k,\R)$ of $G$. Since the elements of $G$ are units in $\A$, this representation has image in $\SL(k,\Z)$. Therefore, $\rho(G)$ is a free abelian subgroups of rank $k-1$ of $\SL(k,\Z)$. Since $\F$ is totally real,  if $g$ is an element of $G\setminus\{1\}$, we see that the norm of $\rho(g^n)$ grows exponentially fast either when $n\to\infty$ or when $n\to -\infty$. So $\rho(g)$ has an eigenvalue of modulus strictly larger than 1.

\end{example}

This example suggests the following open problem which was stated in the arXiv version of \cite{DS1}.

\begin{problem} \rm
Classify compact K\"ahler manifolds of dimension $k\geq 3$ admitting a free abelian group of automorphisms of rank $k-1$ which is of positive entropy.
\end{problem}

In dimension $k=2$, many K3 surfaces and rational surfaces admit positive entropy automorphisms, see e.g. \cite{BK1,BK2, Coble, Diller, Mazur, McMullen1,McMullen2,Uehara} and also \cite{Dolgachev, Keum, Kondo,Oguiso2}. In higher dimension, some partial results on the above problem were obtained in \cite{Zhang4}. We can also ask the same question for groups of rank $k-p$ and for manifolds of dimension large enough, e.g. for groups of rank $k-2$ with $k\geq 4$.

Corollary \ref{c:commute} and the Margulis super-rigidity theorem \cite{Margulis} play a crucial role in a  result of Cantat on the action of a simple lattice on compact K\"ahler manifolds. His result gives an affirmative answer to a version of Zimmer's problem in the case of holomorphic group actions, see also \cite{Zimmer}. The following statement is slightly stronger than the one given in \cite{Cantat2} where the rank bound was $k$ instead of $k-\max(\kappa_X,0)$.

\begin{theorem} \label{t:Cantat}
Let $\Gamma$ be a lattice of a simple algebraic Lie $\R$-group $G$. Assume that $\Gamma$ admits a representation in $\Aut(X)$ with infinite image. Then the real rank of $G$ is at most  $k-\max(\kappa_X,0)$. 
\end{theorem}

We refer to the paper by Cantat \cite{Cantat2} for the proof, see also \cite{CZ2} for more results in this direction.

\medskip\noindent
{\bf Acknowledgments.} I would like to thank Gilles Courtois, Viet-Anh Nguyen, Nessim Sibony, Tuyen Trung Truong for their helps during the preparation of this paper and Benedict Gross for the alternative proof of the property given in Example \ref{ex_torus}. This paper is written for the VIASM Annual Meeting 2012. I would like to thank Ngo Bao Chau, Phung Ho Hai and Le Tuan Hoa for their invitation and hospitality.

\section{Mixed Hodge-Riemann theorem} \label{s:HR}

In this section, we will recall the mixed version of the classical Hodge-Riemann theorem 
which is used for the main results in this paper. We refer to the books by Demailly \cite{Demailly} and Voisin \cite{Voisin} for basic notions and results on Hodge theory.

Let $X$ be a compact K\"ahler manifold of dimension $k$ as above.
For $0\leq p,q\leq k$, denote by $H^{p,q}(X,\C)$ the Hodge cohomology group of bidegree $(p,q)$. We often identify $H^{k,k}(X,\C)$ with $\C$ via the integration of maximal bidegree forms on $X$. 
Define for all $0\leq p\leq k$
$$H^{p,p}(X,\R):=H^{p,p}(X,\C)\cap H^{2p}(X,\R).$$
As a consequence of the Hodge theory, we have 
$$H^{p,p}(X,\C)=H^{p,p}(X,\R)\otimes_\R \C.$$

Denote by $\Kc$ {\it the K\"ahler cone} of $X$, i.e. the set of all classes of K\"ahler forms on $X$, see \cite{DP} for a characterization of $\Kc$. This is a strictly convex open cone in $H^{1,1}(X,\R)$. The closed cone $\overline \Kc$ is called {\it the nef cone}  and its elements are {\it the nef classes} of $X$. 

To each class $\Omega\in H^{k-2,k-2}(X,\R)$, we associate the quadratic form $Q_\Omega$ on $H^{1,1}(X,\R)$ defined by 
 $$Q_\Omega(\alpha,\beta):= - \alpha\smallsmile\beta\smallsmile \Omega \quad \mbox{for}\quad \alpha,\beta\in H^{1,1}(X,\R).$$
For any non-zero class $\Omega'\in H^{k-1,k-1}(X,\R)$, define {\it the primitive subspace} $P_{\Omega'}$ of $H^{1,1}(X,\R)$ associated to $\Omega'$ by 
$$P_{\Omega'}:=\big\{\alpha\in H^{1,1}(X,\R): \ \alpha \smallsmile \Omega' =0 \big\}.$$ 
By Poincar\'e duality, this is a hyperplane of $H^{1,1}(X,\R)$.
We have the following result, see \cite{DN1}. 

\begin{theorem} \label{t:DN}
Let $c_1,\ldots,c_{k-1}$ be K\"ahler classes on $X$. Define $\Omega:=c_1\smallsmile \cdots\smallsmile c_{k-2}$ and 
 $\Omega':=\Omega\smallsmile c_{k-1}$. 
 Let $h^{1,1}$ denote the dimension of $H^{1,1}(X,\R)$.
 Then, the quadratic form $Q_\Omega$ has signature $(h^{1,1}-1,1)$ and is positive definite on $P_{\Omega'}$.
\end{theorem}

When the classes $c_j$ are equal, we obtain the classical Hodge-Riemann theorem. 
The Hodge-Riemann theorem for higher bidegree cohomology groups, 
the hard Lefschetz's theorem and the Lefschetz decomposition theorem can also be generalized in the same way. 

Note that all these results are due to Khovanskii \cite{Khovanskii} and Teissier \cite{Teissier} when $X$ is a projective manifold and the $c_j$'s are integral classes. A linear version was obtained by Timorin in \cite{Timorin}.  
Gromov proved in \cite{Gromov1} that $Q_\Omega$ is semi-positive on $P_{\Omega'}$. His result is in fact enough for our purpose.

\begin{corollary} \label{c:DN}
Let $\alpha$, $\beta$ and $c_j$ be nef classes. Define $\Omega:=c_1\smallsmile \cdots\smallsmile c_{k-2}$. Then, we have
$$|Q_\Omega(\alpha,\alpha)||Q_\Omega(\beta,\beta)|\leq |Q_\Omega(\alpha,\beta)|^2.$$
\end{corollary}
\proof
Note that since $\alpha,\beta$ and the $c_j$'s are nef classes, $Q_\Omega(\alpha,\alpha)$, $Q_\Omega(\beta,\beta)$ and $Q_\Omega(\alpha,\beta)$ are negative. By continuity, we can assume that $\alpha,\beta$ and the $c_j$'s are K\"ahler classes.
Define $A:=|Q(\alpha,\alpha)|$, $B:=|Q(\alpha,\beta)|$, $C:=|Q(\beta,\beta)|$ and $\Omega':=\Omega\smallsmile c_{k-1}$ for some K\"ahler class $c_{k-1}$. 

In order to obtain the corollary, we only have to consider the case where $\alpha$ and $\beta$ are not collinear. So the plane generated by $\alpha$ and $\beta$ intersects $P_{\Omega'}$ along a real line. Let $a\alpha+b\beta$ be a point in this real line with $(a,b)\not=(0,0)$. Since $Q_\Omega$ is positive defined on $P_{\Omega'}$, we deduce that 
$$Aa^2+2Bab+Cb^2=-Q_\Omega(a\alpha+b\beta,a\alpha+b\beta)\leq 0.$$
It follows that $B^2\geq AC$. 
\endproof

Let $\Kc_p$ be the set of all the classes of strictly positive closed $(p,p)$-forms in $H^{p,p}(X,\R)$.
Denote by $\Kc^*_p$ the dual cone of $\Kc_{k-p}$ with respect to the Poincar\'e duality. These cones are strictly convex and open. We also have $\Kc\subset \Kc_1$. 

\begin{definition} \label{d:wHR} \rm
Let $\Theta$ be a class in $\overline\Kc_p^*\setminus\{0\}$ with $p\leq k-2$. 
We say that $\Theta$ is a {\it weak Hodge-Riemann class} (wHR-class for short) if for all K\"ahler classes $c_j$, the quadratic form $Q_\Omega$ is semi-positive on $P_{\Omega'}$, where $\Omega:=\Theta\smallsmile c_1\smallsmile\cdots \smallsmile c_{k-p-2}$ and $\Omega':=\Omega\smallsmile c_{k-p-1}$. 
\end{definition}

By continuity, if $\Theta$ is wHR and the $c_j$'s are nef classes such that $\Omega'\not=0$, then the quadratic form $Q_\Omega$ is still semi-positive on the hyperplane $P_{\Omega'}$. 
Observe also that the set of wHR-classes is closed in $\overline\Kc_p^*\setminus\{0\}$.
By Theorem \ref{t:DN}, if $\Theta$ is a product of nef classes or a limit of such products, then it is a wHR-class. The following result is obtained exactly as in Corollary \ref{c:DN}. 

\begin{proposition}
Let $\Theta$ be a wHR class in $\overline\Kc_p^*\setminus\{0\}$. Let $c_j$ be nef classes and define $\Omega:=\Theta\smallsmile c_1\smallsmile \cdots\smallsmile c_{k-p-2}$. Then, we have 
$$|Q_\Omega(\alpha,\alpha)||Q_\Omega(\beta,\beta)|\leq |Q_\Omega(\alpha,\beta)|^2.$$
\end{proposition}

For any class $\Theta\in \overline\Kc_p^*$, denote by $\overline\Kc(\Theta)$ the closure of $\Theta\smallsmile\Kc$ in $\overline\Kc_{p+1}^*$. We call it {\it the nef cone relative to $\Theta$} or {\it the $\Theta$-nef cone}. 
It is closed, strictly convex and contained in the linear space $\Theta\smallsmile H^{1,1}(X,\R)$. It contains the cone $\Theta\smallsmile \overline \Kc$ and can be strictly larger than this  cone; in other words, $\Theta\smallsmile \overline\Kc$ is not always closed. Observe that if $\Theta$ is a wHR-class, so are the classes in $\overline\Kc(\Theta)\setminus\{0\}$.

\begin{proposition} \label{p:wHR_set}
Let $\pi:X'\to X$ be a holomorphic map between compact K\"ahler manifolds. If $\Theta'$ is a wHR-class on $X'$ such that $\pi_*(\Theta')\not=0$, then $\pi_*(\Theta')$ is a wHR-class on $X$. In particular, if $\Theta$ is the class of an irreducible analytic subset of $X$ then $\Theta$ is wHR.
\end{proposition}
\proof
If $c$ is a nef class on $X$ then $\pi^*(c)$ is a nef class on $X'$. Therefore, the first assertion is a direct consequence of Definition \ref{d:wHR}. For the second assertion, assume that $\Theta$ is the class of an irreducible analytic subset of $X$. If the analytic set is smooth, it is enough to apply Theorem \ref{t:DN} to this manifold. Otherwise, we use a resolution of singularities and the first assertion allows us to reduce the problem to the smooth case.
\endproof

\begin{lemma} \label{l:nef_wedge}
Let $\Theta$ be a class in $\overline\Kc_p^*$, $M$ a class in $H^{q,q}(X,\R)$ and 
$\Theta'$ a class in $\overline\Kc(\Theta)$. Assume that $\Theta\smallsmile M$ is in $\overline\Kc_{p+q}^*$ and write 
 $\Theta'=\Theta\smallsmile L$ with $L\in H^{1,1}(X,\R)$. Then $\Theta\smallsmile M\smallsmile L$ 
 is a class in  $\overline\Kc(\Theta\smallsmile M)$ and it
 does not depend on the choice of $L$.  
\end{lemma}
\proof
The independence of the choice of $L$ is obvious. Since $\Theta'$ is in $\overline\Kc(\Theta)$ there is a sequence of K\"ahler classes $L_n$ such that $\Theta\smallsmile L_n\rightarrow \Theta\smallsmile L$. We deduce that $\Theta\smallsmile M\smallsmile L_n$ converge to $\Theta\smallsmile M\smallsmile L$. Therefore, the last class belongs to $\overline \Kc(\Theta\smallsmile M)$.
\endproof

\begin{definition} \rm 
Let $\Theta$ and $\Theta'$ be two classes in $H^{p,p}(X,\R)$. We say that they are {\it numerically almost equivalent} and we write $\Theta\simeq_n\Theta'$ if 
$$(\Theta-\Theta')\smallsmile c_1\smallsmile \ldots\smallsmile c_{k-p}=0$$
for all classes $c_j\in H^{1,1}(X,\R)$. 
\end{definition} 

Note that for $\Theta$ in $\overline\Kc_p^*$ we have $\Theta\not\simeq_n 0$ if and only if $\Theta\not=0$. We will need the following proposition which uses the ideas from \cite{DS1,DS3,Zhang3}. 

\begin{proposition} \label{p:HR}
Let $\Theta$ be a wHR-class in $\overline\Kc_p^*$.
Let $\Theta_1$ and $\Theta_2$ be two classes in $\overline \Kc(\Theta)$. 
Write $\Theta_j=\Theta\smallsmile L_j$ with $L_j\in H^{1,1}(X,\R)$.
Assume that  $\Theta\smallsmile L_1\smallsmile L_2 = 0$. Then, we also have $\Theta\smallsmile L_1^2=\Theta\smallsmile L_2^2=0$. Moreover, there is a pair of  real numbers $(t_1,t_2)\not=(0,0)$ such that 
$ \Theta\smallsmile (t_1L_1 + t_2 L_2) \simeq_n 0.$
\end{proposition}
\proof
By Lemma \ref{l:nef_wedge}, the classes 
$\Theta\smallsmile L_1\smallsmile L_2$ and $\Theta\smallsmile L_j^2$ belong to $\overline\Kc^*_{p+2}$ and depend only on $\Theta_j$ but not on the choice of $L_j$. 
Observe also that we only need to consider the case where $\Theta_1$ and $\Theta_2$ are linearly independent. So  $\Theta_1,\Theta_2$ belong to $\overline \Kc_{p+1}^*\setminus\{0\}$ and $L_1,L_2$ are linearly independent. 
Denote by $H$ the real plane in $H^{1,1}(X,\R)$ generated by $L_1$ and $L_2$. 

Let $c_1,\ldots, c_{k-p-1}$ be K\"ahler classes. 
Define $\Omega:=\Theta\smallsmile c_1\smallsmile \cdots\smallsmile c_{k-p-2}$ and $\Omega':=\Omega\smallsmile c_{k-p-1}$. We deduce from the hypothesis on
$\Theta\smallsmile L_1\smallsmile L_2$ and the definition of $Q_\Omega$ that $Q_\Omega$ is semi-negative on $H$. Since $\Theta$ is a  wHR-class, $Q_\Omega$ is semi-positive on $P_{\Omega'}$. It follows that $Q_\Omega$ vanishes on $H\cap P_{\Omega'}$. 

Consider a pair of real numbers $(t_1,t_2)\not=(0,0)$ such that $t_1L_1+t_2L_2$ belongs to $P_{\Omega'}$. We have
$$\Theta \smallsmile (t_1L_1+t_2L_2)^2\smallsmile c_1\smallsmile \ldots\smallsmile c_{k-p-2}=0.$$
Hence,
$$\Theta \smallsmile (t_1^2L_1^2+t_2^2L_2^2)\smallsmile c_1\smallsmile\ldots\smallsmile c_{k-p-2}=0.$$
Since $\Theta\smallsmile L_j^2\in\overline\Kc_{p+2}^*$, we conclude that $\Theta\smallsmile L_j^2=0$ if $t_j\not=0$. 

Recall that we suppose $\Theta_j\in \overline \Kc_{p+1}^*\setminus\{0\}$. Therefore, $L_j$ does not belong to $P_{\Omega'}$. Thus, $t_j\not=0$ and  $\Theta\smallsmile L_j^2=0$. We can assume that $t_1=1$ and write $t:=t_2$. The number $t$ is the unique real number such that $L_1+tL_2\in H\cap P_{\Omega'}$. By Cauchy-Schwarz's inequality applied to the restriction of $Q_\Omega$ to $P_{\Omega'}$, we have $Q_\Omega(L_1+tL_2,c)=0$ for every $c$ in the hyperplane $P_{\Omega'}$. This together with the first assertion in the proposition implies that $Q_\Omega(L_1+tL_2,c)=0$ for every $c\in H^{1,1}(X,\R)$. Then, it follows from the Poincar\'e duality that 
$$\Theta\smallsmile (L_1+tL_2)\smallsmile c_1\smallsmile \ldots\smallsmile c_{k-p-2}=0.$$

In order to obtain the last assertion in the proposition, it remains to check that $t$ is independent of $c_j$. 
By definition, $t$ depends symmetrically on the $c_j$'s. However, the last identity, which is stronger than the property $L_1+tL_2\in P_{\Omega'}$, shows that it does not depend on $c_{k-p-1}$.  We conclude that $t$ does not depend on $c_j$ for every $j$. This completes the proof of the proposition.
\endproof

\section{Topological entropy and dynamical degrees} \label{s:entropy}

The topological entropy of a map is a dynamical invariant. In the case of holomorphic automorphisms of a compact K\"ahler manifold, results by Gromov and Yomdin imply that the topological entropy is in fact an algebraic invariant.

Let $f$ be a holomorphic automorphism of $X$. It defines a dynamical system on $X$.  
Denote by $f^n:=f\circ\cdots\circ f$, $n$ times, the iterate of order $n$ of $f$. 
If $x$ is a point of $X$, the orbit of $x$ is the sequence of points 
$$x,f(x),f^2(x),\ldots, f^n(x),\ldots$$
The topological entropy of $f$ measures the divergence of the orbits or in some sense it measures the rate of expansion of the "number" of orbits one can distinguish when the time $n$ goes to infinity. The formal definition is given below.

\begin{definition} \label{d:separate} \rm
Let $\epsilon>0$ and $n\in\N$. Two points $x$ and $y$ in $X$ are said to be {\it $(n,\epsilon)$-separated} if we have 
 for some integer $0\leq j\leq n-1$
 $$\dist(f^j(x),f^j(y))>\epsilon.$$ 
\end{definition}

\begin{definition}  \label{d:entropy} \rm
Let $N(\epsilon,n)$ denote the maximal number of points mutually $(n,\epsilon)$-separated. {\it The (topological) entropy} of $f$ is given by the formula
$$h_t(f):=\sup_{\epsilon>0} \limsup_{n\to\infty} {\log N(\epsilon,n) \over n}=\lim_{\epsilon\to 0} \limsup_{n\to\infty} {\log N(\epsilon,n) \over n}\cdot$$
\end{definition}

Note that the notion  of separated points depends on the metric on $X$ but the topological entropy is independent of the choice of the metric. So the topological entropy is a topological invariant.

The pull-back operator $f^*$ on differential forms induces a graded automorphism of the Hodge cohomology ring of $X$
$$f^*:\bigoplus_{0\leq p,q\leq k} H^{p,q}(X,\C)\to \bigoplus_{0\leq p,q\leq k} H^{p,q}(X,\C).$$
A similar property holds for the de Rham cohomology ring which can be identified to the real part of the Hodge cohomology ring. We have a graded automorphism
$$f^*:\bigoplus_{0\leq m \leq 2k} H^m(X,\R)\to \bigoplus_{0\leq m\leq 2k} H^m(X,\R).$$
The last operator preserves a lattice on which it is identified to
$$f^*:\bigoplus_{0\leq m \leq 2k} {H^m(X,\Z) \over \mbox{torsion}}\to \bigoplus_{0\leq m\leq 2k} {H^m(X,\Z)\over \mbox{torsion}}\cdot$$
In particular, on a suitable basis, the map $f^*:H^m(X,\R)\to H^m(X,\R)$ is given by a square matrix with integer entries. When $X$ is a projective manifold, we can also consider the action of $f$ on N\'eron-Severi groups.

The following algebraic invariants were implicitly considered in Gromov \cite{Gromov2}, see also \cite{DS2, Friedland, RS}.

\begin{definition} \label{d:dyn_degree}\rm
We call {\it dynamical degree of order $p$ of $f$} the spectral radius $d_p(f)$ of the linear morphism $f^*:H^{p,p}(X,\C)\to H^{p,p}(X,\C)$ and we call {\it algebraic entropy of $f$} the number 
$$h_a(f):=\max_{0\leq p\leq k} \log d_p(f).$$
\end{definition}

It follows from the above discussion that each $d_p(f)$ is a root of a monic polynomial with integer coefficients. In particular, it is an algebraic number.

If $\omega$ is a K\"ahler form on $X$, it is not difficult to see that the dynamical degrees can by computed with the formula
\begin{eqnarray} \label{e:degree}
d_p(f) & = & \lim_{n\to\infty} \Big[\int_X (f^n)^*(\omega^p)\wedge \omega^{k-p}\Big]^{1/n}
=\lim_{n\to\infty} \Big[(f^n)^*\{\omega\}^p\smallsmile \{\omega\}^{k-p}\Big]^{1/n} \nonumber\\
& = & \lim_{n\to\infty} \Big[\int_X \omega^p\wedge (f^n)_*(\omega^{k-p})\Big]^{1/n}
=\lim_{n\to\infty} \Big[\{\omega\}^p\smallsmile (f^n)_*\{\omega\}^{k-p}\Big]^{1/n}.
\end{eqnarray}
We see that  $d_0(f)=d_k(f)=1$ and $d_p(f)=d_{k-p}(f^{-1})$.
 
The following result is a consequence of a theorem by Gromov \cite{Gromov2} and another theorem by Yomdin \cite{Yomdin}. It shows that the topological entropy of a holomorphic automorphism can be computed algebraically.

\begin{theorem} \label{t:GY}
We have
$$h_t(f)=h_a(f).$$
\end{theorem}

Gromov theorem implies that $h_t(f)\leq h_a(f)$. A similar property holds for all meromorphic self-maps on compact K\"ahler manifolds, see \cite{DS2}. 

Yomdin theorem says that if $V\subset X$ is a real manifold, smooth up to the boundary, then the volume growth of the sequence $f^n(V)$, $n\geq 0$, is bounded by $h_t(f)$. More precisely, we have
$$h_t(f)\geq \limsup_{n\to\infty} {1\over n}\log \volume f^n(V),$$
where we use  the $m$-dimensional volume $\volume(\cdot)$ with $m:=\dim_\R V$. 
Yomdin's theorem holds for all smooth maps on compact real manifolds. 

Applying Yomdin theorem to real compact manifolds without boundary in $X$, we obtain that 
$$h_t(f)\geq\log \rho_m(f)$$ 
if $\rho_m(f)$ is the spectral radius of $f^*:H^m(X,\R)\to H^m(X,\R)$. In particular, we obtain the reverse of the above Gromov's inequality. 

Let $\rho_{p,q}(f)$ denote the spectral radius of $f^*:H^{p,q}(X,\C)\to H^{p,q}(X,\C)$. Arguing as above, we obtain
$$h_t(f)\geq \log \rho_{p,q}(f).$$
This together with Gromov's inequality yields
$$\rho_{p,q}(f)\leq \max_{0\leq p\leq k} d_p(f).$$
In fact, the following more general result holds, see \cite{Dinh}.

\begin{proposition} \label{p:deg_p_q}
We have for $0\leq p,q\leq k$
$$\rho_{p,q}(f)\leq \sqrt{d_p(f)d_q(f)}.$$
\end{proposition}

This inequality also explains why the dynamical degrees $d_p(f)$ play a more important role than the degrees $\rho_{p,q}(f)$, $p\not=q$, in the dynamical study of $f$. A similar property holds for general meromorphic self-maps on compact K\"ahler manifolds. 

The following result is a direct consequence of Corollary \ref{c:DN} and the identity (\ref{e:degree}).

\begin{proposition} \label{p:deg_concave}
The function $p\mapsto \log d_p(f)$ is concave in $p$, that is, 
$$d_p(f)^2\geq d_{p-1}(f)d_{p+1}(f)\quad \mbox{for} \quad 1\leq p\leq k-1.$$
In particular, there are two integers $0\leq r\leq s\leq k$ such that 
$$1=d_0(f)<\cdots <d_r(f)=\cdots=d_s(f)>\cdots>d_k(f)=1.$$ 
\end{proposition}

We obtain the following corollary, see \cite{DS1}.

\begin{corollary} \label{c:entropy_bound}
The automorphism $f$ has positive entropy if and only if $d_1(f)>1$ (resp. $d_{k-1}(f)>1$). Moreover, in this case, there is a number $A>1$ depending only on the second Betti number of $X$  such that
$$h_t(f)\geq \log A  \quad \mbox{and} \quad d_p(f)\geq A \quad \mbox{for} \quad 1\leq p\leq k-1.$$
\end{corollary}
\proof
The first assertion is a consequence of Theorem \ref{t:GY} and Proposition \ref{p:deg_concave}. Assume now that $f$ has positive entropy.
We have $d_p(f)>1$ for $1\leq p\leq k-1$. 
It suffices to prove that $d_1(f)\geq A$ and  $d_{k-1}(f)\geq A$. We only have to check the first inequality since $d_{k-1}(f)=d_1(f^{-1})$.  

It follows from Proposition \ref{p:deg_p_q} that $d_1(f)$ is the spectral radius of $f^*:H^2(X,\R)\to H^2(X,\R)$. So it is the largest root of a monic polynomial of degree $b_2$ with integer coefficients, where $b_2:=\dim H^2(X,\R)$ denotes the second Betti number of $X$. 
If this polynomial admits a coefficient of absolute value larger than $2^{b_2}b_2!$ then it admits a root of modulus larger than 2. In this case, we have $d_1(f)\geq 2$. Otherwise, the polynomial belongs to a finite family and hence $d_1(f)$ belongs to a finite set depending only on $b_2$.
The result follows.
\endproof

\section{Fibrations and relative dynamical degrees} \label{s:fibration}

We will consider the restriction of an automorphism to analytic sets which may be singular. In general, a resolution of singularities gives us maps which are no more holomorphic. So it is useful to extend the notion of dynamical degrees to meromorphic maps. 

For the moment, let $(X,\omega)$ be a compact K\"ahler manifold. Let $f:X\to X$ be a meromorphic map which is dominant, i.e. its image contains a Zariski open subset of $X$. {\it The dynamical degree of order $p$} of $f$ is defined by
 \begin{equation} \label{e:degree_bis}
d_p(f)=\lim_{n\to\infty} \Big[\int_X (f^n)^*(\omega^p)\wedge \omega^{k-p}\Big]^{1/n}.
\end{equation}

It is not difficult to see that the definition does not depend on the choice of $\omega$. The existence of the above limit is not obvious. It is based on some result on the regularization of positive closed currents \cite{DS2}, see also \cite{Demailly2}. Dynamical degrees are bi-meromorphic invariants. More precisely, we have the following result, see \cite{DS2}.

\begin{theorem} \label{t:deg_bimer}
Let $f$ and $g$ be dominant meromorphic self-maps on compact K\"ahler manifolds $X$ and $Y$ respectively, of the same dimension $k$. Let $\pi:X\to Y$ be a bi-meromorphic map. Assume that $g\circ \pi=\pi\circ f$. Then, we have $d_p(f)=d_p(g)$ for $0\leq p\leq k$. 
\end{theorem}

So we can extend the notion of dynamical degrees to meromorphic maps on varieties by using a resolution of singularities. 
Note that Proposition \ref{p:deg_concave} still holds in this case except we only have $d_k(f)\geq 1$ with equality when $f$ is a bi-meromorphic map.

The last theorem can be viewed also as a consequence of Theorem \ref{t:deg_fibration} below. Consider now a dominant meromorphic map $g:Y\to Y$, where $Y$ is a compact K\"ahler manifold of dimension $l\leq k$. Let $\pi:X\to Y$ a dominant meromorphic map and assume as above that  $g\circ \pi=\pi\circ f$. So $f$ preserves the meromorphic fibration defined by $\pi$. 

We can define {\it the dynamical degree of $f$ relative to the fibration} by
$$d_p(f|\pi):=\lim_{n\to\infty} \Big[\int_{\pi^{-1}(y)} (f^n)^*(\omega^p)\wedge \omega^{k-l-p}\Big]^{1/n},$$
where $y$ is a generic point in $Y$. The definition does not depend on the generic choice of  $y$ and the function $p\mapsto \log d_p(f|\pi)$ is concave on $p$.
The following result relates the dynamical degrees of $f$ and the ones of $g$, see \cite{DNT}.

\begin{theorem} \label{t:deg_fibration}
Let $f,g,\pi$ be as above. Then, we have for $0\leq p\leq k$
$$d_p(f)=\max_{\max(0,p-k+l)\leq s\leq \min(p,l)} d_s(g)d_{p-s}(f|\pi).$$
\end{theorem}

Note that the domain of $s$ in the last formula is exactly the set of $s$ such that $d_s(g)$ and $d_{p-s}(f|\pi)$ are meaningful. In the case where $k=l$, we necessarily have $s=p$ and $d_0(f|\pi)=1$. So the last formula implies Theorem \ref{t:deg_bimer}. We will  apply the last theorem to the case of pluricanonical fibrations of $X$. 

Let $K_X$ denote the canonical line bundle of $X$. Let $H^0(X,K_X^N)$
denote the space of holomorphic sections of $K_X^N$ and 
$H^0(X,K_X^N)^*$ its dual space. Assume that $H^0(X,K_X^N)$ has a
positive dimension. If $x$ is a generic point in $X$, the family $H_x$ of
sections which vanish at $x$ is a hyperplane of  $H^0(X,K_X^N)$
passing through 0. So
the correspondence $x\mapsto H_x$ defines 
a meromorphic map
$$\pi_N:X\rightarrow \P H^0(X,K_X^N)^*$$
from $X$ to the projectivization of $H^0(X,K_X^N)^*$ which is called a
{\it pluricanonical fibration} of $X$. Let $Y_N$ denote the image of $X$ by
$\pi_N$. The {\it Kodaira dimension} of $X$ is
$\kappa_X:=\max_{N\geq 1}\dim Y_N$. When $H^0(X,K_X^N)=0$ for
every $N\geq 1$, the Kodaira dimension of $X$ is defined to be $-\infty$. 
We have the following result, see \cite{NZ, Ueno}.

\begin{theorem} \label{t:Iitaka_fix}
Let $f:X\rightarrow X$ be a dominant meromorphic map. Assume that
$\kappa_X\geq 1$. Then $f$ preserves the pluricanonical fibration
$\pi_N:X\rightarrow Y_N$. Moreover, the map $g_N:Y_N\rightarrow Y_N$
induced by $f$ is periodic, i.e. $g_N^m=\id$ for some integer $m\geq
1$. 
\end{theorem}

We deduce that $d_p(g_N)=1$ for every $p$. This property can also be deduced from a weaker property that $g_N$ is the restriction to $Y_N$ of a linear map on $\P H^0(X,K_X^N)^*$ which is a consequence of the definition of $g_N$.
The following result is a consequence of Corollary \ref{c:entropy_bound} and Theorem \ref{t:deg_fibration}.

\begin{corollary} \label{c:deg_Iitaka}
Let $f$ be a holomorphic automorphism of $X$. Assume that $0\leq \kappa_X\leq k-1$. Let $Y_N,\pi_N,g_N$ be as above. Then 
$$h_t(f)=\max_{1\leq p\leq k-\dim Y_N-1} d_p(f|\pi_N) \quad \mbox{and} \quad d_1(f)=d_1(f|\pi_N).$$
In particular, $f$ has positive entropy if and only if $d_1(f|\pi_N)\geq A$, where $A>1$ is the constant given in Corollary \ref{c:entropy_bound}.  
\end{corollary}

\section{Tits alternative for automorphism groups} \label{s:Tits}

We are now ready to give the proof of Theorem \ref{t:Zhang}. 
We first recall two important results due to Fujiki and Lieberman \cite{Fujiki,Lieberman}. 

\begin{theorem} \label{t:FL1}
Let $c$ be a K\"ahler class on $X$. Denote by $\Aut_c(X)$ the group of elements $g$ of $\Aut(X)$ such that $g^*(c)=c$. Then $\Aut_c(X)$ is a finite union of connected components of $\Aut(X)$. 
\end{theorem}

Let $\Alb(X)$ denote the Albanese torus of $X$ and $\phi:X\to \Alb(X)$ the Albanese map. The identity component of the automorphism group of $\Alb(X)$ is denoted by $A(X)$. This is the group of translations on $\Alb(X)$ which is isomorphic to $\Alb(X)$ as complex Lie groups. It is not difficult to see that any automorphism $g$ of $X$ induces an automorphism $h$ of $\Alb(X)$ such that $h\circ\phi=\phi\circ g$. So we have a natural Lie group morphism 
$$\psi:\Aut_0(X)\to A(X).$$ 

\begin{theorem} \label{t:FL2}
The kernel $\ker(\psi)$ of $\psi$ is a linear algebraic $\C$-group. 
\end{theorem}

We recall also the following version of Tits' theorem \cite{Tits}.

\begin{theorem} \label{t:Tits}
Let $G$ be a linear $\R$-group. Then, it satisfies the Tits alternative, that is, any subgroup of $G$ either has a free non-abelian subgroup or virtually solvable, i.e. possesses a solvable subgroup of finite index. 
\end{theorem}

The following result gives us the first assertion of Theorem \ref{t:Zhang}. It was obtained in \cite{CWZ,Cantat}.

\begin{theorem} \label{t:Tits_Aut}
The group $\Aut(X)$ satisfies the Tits alternative.
\end{theorem}

We first prove a preliminary lemma.

\begin{lemma} \label{l:vir_soluble}
Let $A$ be a group and $B$ a normal subgroup of $A$. If $B$ and $A/B$ are virtually solvable, then $A$ is virtually solvable. If $B$ and $A/B$ satisfy the Tits alternative, then $A$ satisfies also the Tits alternative.
\end{lemma}
\proof
The second assertion is a direct consequence of the first one. 
Assume that $B$ and $A/B$ are virtually solvable. We show that $A$ is virtually solvable.
Let $\pi:A\to A/B$ be the canonical group morphism. If $D$ is a solvable finite index subgroup of $A/B$, we can replace $A$ by $\pi^{-1}(D)$ in order to assume that $A/B$ is solvable.  Let 
$$\{1\}=D_0\triangleleft D_1 \triangleleft \cdots \triangleleft D_{m-1}\triangleleft D_m= A/B$$
be a subnormal series such that $D_j$ is normal in $A/B$ and $D_{j+1}/D_j$ is abelian for every $0\leq j\leq m-1$. We can use here the derived series of $A/B$.  

By lattice theorem (correspondence theorem), there is a subnormal series
$$B=B_0\triangleleft B_1 \triangleleft \cdots \triangleleft B_{m-1}\triangleleft B_m= A$$
such that $B_j$ is normal in $A$ and $B_j/B_{j-1}=D_j/D_{j-1}$ for $1\leq j\leq m$. 
Recall that $B$ is virtually solvable. We will show that $B_1$ satisfies the same property and then using a simple induction, we obtain that $A$ is virtually solvable. So in order to simplify the notation, we can assume that $m=1$ or equivalently $A/B$ is abelian.  

Let $C$ be a solvable finite index subgroup of $B$. We can replace $C$ by the intersection of $bCb^{-1}$ with $b\in B$ in order to assume that $C$ is normal in $B$. Without loss of generality, we can also assume that $C$ is a maximal normal solvable subgroup of $B$ with finite index. 
The maximality and the lattice theorem imply that $B/C$ admits no solvable normal subgroup different from $\{1\}$. 
Since $B$ is normal in $A$, we have $a^{-1}Ca\subset B$ for every $a\in A$. We claim that $C$ is normal in $A$, i.e. $a^{-1}Ca=C$ for every $a\in A$.

Taking into account this property, we first complete the proof of the lemma. 
Observe that if $a$ is an element of $A$ then $b\mapsto a^{-1}ba$ induces an automorphism of the group $B/C$.
Let $A'$ denote the set of all elements $a\in A$ such that the above automorphism is identity. 
Since $B/C$ is a finite group, $A'$ is a finite index subgroup of $A$. 

Since $A/B$ is abelian, $A'':=[A',A']$ is a subgroup of $B$. By construction, if $a'$ is an element of $A'$ and $b$ an element of $B$, then $[a',b]$ is an element of $C$. We deduce that $[A'',A'']$ is a subgroup of $C$; in particular, it is solvable. Thus, $A'$ is solvable. It remains to prove the above claim. 

Define $D:=a^{-1}Ca$. Since $b\mapsto a^{-1}ba$ is an automorphism of $B$, $D$ is a maximal normal  solvable subgroup of $B$ with finite index and $B/D$ is isomorphic to $B/C$. So it suffices to check that $D\subset C$. The natural short exact sequence 
$$\{1\}\longrightarrow D \longrightarrow B\longrightarrow B/D\longrightarrow \{1\}$$
induces the following one
$$\{1\}\longrightarrow {D\over C\cap D} \stackrel{\pi_1}{\longrightarrow} {B\over C\cap D} \longrightarrow {B\over D} \longrightarrow \{1\}.$$
Similarly, we have 
$$\{1\}\longrightarrow {C\over C\cap D} \longrightarrow {B\over C\cap D} \stackrel{\pi_2}{\longrightarrow} {B\over C} \longrightarrow \{1\}.$$

Since $\pi_1$ is injective and $\pi_2$ is surjective, the image of $\pi_2\circ\pi_1$ is a normal subgroup of $B/C$. On the other hand, this subgroup should be solvable since $D$ is solvable. We deduce from the maximality of $C$ that the image of $\pi_2\circ\pi_1$ is equal to $\{1\}$. Hence, $D\subset C$. This completes the proof of the lemma.
\endproof

\noindent
{\bf Proof of Theorem \ref{t:Tits_Aut}.} 
Let $G$ be a subgroup of $\Aut(X)$ which does not contain any free non-abelian subgroup. We have to show that $G$ admits a solvable subgroup of finite index.  
Consider the natural group morphism
$$\rho:\Aut(X)\to \GL(H^2(X,\R)).$$
By Theorem \ref{t:Tits}, 
$\rho(G)$ is virtually solvable. By Lemma \ref{l:vir_soluble}, we only have to check that $G\cap\ker\rho$ is virtually solvable.  

By Theorem \ref{t:FL1}, $G\cap \ker(\rho)$ is a finite extension of $G\cap\Aut_0(X)$. So we only have to check that $G\cap\Aut_0(X)$ is virtually solvable.  
Since $\psi(G\cap\Aut_0(X))$ is abelian, by Lemma \ref{l:vir_soluble}, it suffices to show that $\ker\psi\cap G$ is virtually solvable. But this is a consequence of Theorems \ref{t:FL2} and \ref{t:Tits}. 
\hfill $\square$ \noindent

\bigskip

We now turn to the proof of the second assertion in Theorem \ref{t:Zhang}. 
Let $G$ be a group as in this theorem.
By Theorem \ref{t:Tits_Aut}, $G$ is virtually solvable. We will need the following version 
of the Lie-Kolchin theorem due to Keum-Oguiso-Zhang \cite{KOZ}.

\begin{theorem} \label{t:LK}
Let $H$ be a virtually solvable group acting linearly on a strictly convex closed cone $C$ of finite dimension. Then, $H$ admits a finite index subgroup $H'$  and a non-zero vector $v\in C$ such that the half-line $\R_+v$ is invariant by $H'$. 
\end{theorem} 

This result was obtained by induction on the derived length of a suitable finite index solvable subgroup of $G$. The case where $G$ is abelian is a version of the classical Perron-Frobenius theorem.

\medskip

Observe that in the case $\kappa_X=k$ the second assertion in Theorem \ref{t:Zhang} is a direct consequence of Theorem \ref{t:Iitaka_fix} since in this case every automorphism has zero entropy. Assume now that $\kappa_X\leq k-1$. 
Fix now an integer $N$ such that $\dim Y_N=\kappa_X$, where $Y_N$ is as defined in Section \ref{s:fibration}.
In order to simplify the notation, define $\pi:=\pi_N$, $Y:=Y_N$ and $\kappa:=\max(\kappa_X,0)$. If $\kappa_X=-\infty$, we consider that $Y$ is a point. Let $\Theta_\kappa$ denote the class of a generic fiber of $\pi$. 
In general, the generic fibers of $\pi$ are not necessarily irreducible. However, by Stein's factorization theorem \cite[Ch. 10.6]{GR}, their irreducible components have the same cohomology class.
Therefore,
by Proposition \ref{p:wHR_set},
$\Theta_\kappa$ is a wHR-class in $\overline\Kc^*_\kappa\setminus\{0\}$. 
By Theorem \ref{t:Iitaka_fix},
this class is fixed under the action of $\Aut(X)$. 

\begin{lemma} \label{l:inv_seq}
There is a finite index subgroup $G'$ of $G$ such that 
for every $\kappa\leq p\leq k-1$, there exists a wHR-class $\Theta_p$ in $\overline\Kc^*_p\setminus\{0\}$ and a character $\chi_p:G'\to \R^*$ of $G'$ such that $g^*(\Theta_p)=\chi_p(g)\Theta_p$ for $g\in G'$. Moreover, we have $\Theta_p\in \overline\Kc(\Theta_{p-1})$ when $p\geq \kappa+1$. 
\end{lemma}
\proof
We construct  $\Theta_p$ by induction on $p$. 
The class $\Theta_\kappa$ was already constructed above and we can take $G'=G$. Assume that $\Theta_{p-1}$ was constructed. Then,  $G'$ induces a linear action on the strictly convex cone $\overline\Kc(\Theta_{p-1})$. 
By Theorem \ref{t:LK}, replacing $G'$ by a suitable finite index subgroup, we can find a class $\Theta_p\in \overline\Kc(\Theta_{p-1})\setminus\{0\}$ whose direction is invariant by $G'$. Since $\Theta_{p-1}$ is a wHR-class,  $\Theta_p$ is also a wHR-class. 
\endproof

Consider the group morphism $\phi:G'\to \R^{k-\kappa-1}$ given by
$$\phi(g):=\big(\log\chi_{\kappa+1}(g),\ldots, \log\chi_{k-1}(g)\big).$$
The following lemma will permit to show that $\Im(\phi)$ is discrete. 

\begin{lemma} \label{l:phi_bound}
We have $\|\phi(g)\|\geq {1\over 2}\log d_{k-1}(g)$ for all $g\in G'$. 
\end{lemma}
\proof
Assume that $\|\phi(g)\|< {1\over 2}\log d_{k-1}(g)$ for some $g\in G'$. 
Then, we have 
\begin{equation} \label{e:phi_bound}
d_{k-1}(g)^{-1/2}<\chi_p(g)<d_{k-1}(g)^{1/2}
\end{equation}
for every $p$. 
Recall that $d_{k-1}(g)=d_1(g^{-1})=d_1(g^{-1}|\pi)$, see Corollary \ref{c:deg_Iitaka}. Let $\Theta_p$ be as in Lemma \ref{l:inv_seq} and write $\Theta_p=\Theta_{p-1}\smallsmile L_p$ with some class $L_p\in H^{1,1}(X,\R)$. 

Since $g^{-1}$ preserves $\overline\Kc(\Theta_\kappa)$, it follows from the  Perron-Frobenius theorem that there is a class $\Theta \in\overline\Kc(\Theta_\kappa)\setminus\{0\}$ depending on $g$ such that 
$$(g^{-1})^*(\Theta)=d_1(g^{-1}|\pi)\Theta=d_{k-1}(g)\Theta$$
or equivalently
 $$g^*(\Theta)=d_{k-1}(g)^{-1}\Theta.$$
Write $\Theta=\Theta_\kappa\smallsmile L$ with $L\in H^{1,1}(X,\R)$.

By Lemma \ref{l:nef_wedge}, $\Theta_p\smallsmile L$ does not depend on the choice of $L$ and  it is not difficult to see that 
\begin{equation} \label{e:smile_L}
g^*(\Theta_p\smallsmile L)=\chi_p(g)d_{k-1}(g)^{-1}\Theta_p\smallsmile L.
\end{equation}
Since $g^*=\id$ on $H^{k,k}(X,\R)$ and $\chi_{k-1}(g)d_{k-1}(g)^{-1}\not=1$, we deduce that $\Theta_{k-1}\smallsmile L=0$. Let $q\leq k-1$ be the smallest integer such that $\Theta_q\smallsmile L=0$.  

Since $\Theta$ belongs to $\overline\Kc^*_{\kappa+1}\setminus \{0\}$, we have $\Theta\not\simeq_n 0$. Therefore, we have $q\geq \kappa+1$. 
We  have $\Theta_{q-1}\smallsmile L_q\smallsmile L= 0$. By Proposition \ref{p:HR}, there is a pair of real numbers $(t_1,t_2)\not=(0,0)$ such that 
$$\Theta_{q-1}\smallsmile (t_1L_q+t_2L)\simeq_n 0.$$ 
Using the action of $g^*$ and the relation (\ref{e:smile_L}), we obtain that 
$$\Theta_{q-1}\smallsmile \big(t_1\chi_q(g)L_q+t_2\chi_{q-1}(g)d_{k-1}(g)^{-1}\big)\simeq_n0.$$
The last two identities together with (\ref{e:phi_bound}) yield 
$\Theta_{q-1}\smallsmile L\simeq_n 0$. By Lemma \ref{l:nef_wedge}, $\Theta_{q-1}\smallsmile L$ belongs to $\overline\Kc_q^*$. Thus,
$\Theta_{q-1}\smallsmile L=0$. This contradicts the minimality of $q$. 
The lemma follows.
\endproof

\noindent
{\bf End of the proof of Theorem \ref{t:Zhang}.}
When $h_t(g)=0$, the spectral radius of $g^*$ on $\oplus H^m(X,\R)$ is equal to 1. Since $g^*$ is given by a matrix with integer entries, we deduce that all eigenvalues of $g^*$ have modulus 1. It follows that $\phi(g)=0$. Conversely, if $\phi(g)=0$, by Lemma \ref{l:phi_bound}, $d_{k-1}(g)=1$; thus, by Corollary \ref{c:entropy_bound}, we get 
$h_t(g)=0$. So we have $N'=\ker \phi$. In particular, $N'$ is a normal subgroup of $G'$. The group $G'/N'$ is isomorphic to $\phi(G')$. 
By Corollary \ref{c:entropy_bound} and Lemma \ref{l:phi_bound}, $\phi(G')$ is a discrete subset of $\R^{k-\kappa-1}$. So $G'/N'$ is a free abelian group of rank $\leq k-\kappa-1$. This finishes the proof of the theorem.
\hfill $\square$

\noindent
T.-C. Dinh, UPMC Univ Paris 06, UMR 7586, Institut de
Math{\'e}matiques de Jussieu, 4 place Jussieu, F-75005 Paris, France.\\ 
DMA, UMR 8553, Ecole Normale Sup\'erieure,
45 rue d'Ulm, 75005 Paris, France.\\
{\tt dinh@math.jussieu.fr}, {\tt http://www.math.jussieu.fr/$\sim$dinh}

 \end{document}